\documentclass[letterpaper, 10 pt, conference]{ieeeconf}

\IEEEoverridecommandlockouts 
\usepackage{amsmath, amssymb, mathtools, mathrsfs, xcolor}
\usepackage[american]{circuitikz}
\usepackage{comment}
\usepackage{subcaption}
\usepackage{pgfplots}
\usepackage{tikz}

\ctikzset{bipoles/length=1.4cm}

\newcommand{\bbR}{\mathbb{R}}

\newcommand{\bbZ}{\mathbb{Z}}

\newcommand{\calY}{\mathcal{Y}}

\newcommand{\calU}{\mathcal{U}}

\renewcommand{\epsilon}{ \varepsilon}

\newcommand{\bbm}{\begin{bmatrix*}}
\newcommand{\ebm}{\end{bmatrix*}}

\newcommand{\abs}[1]{\lvert#1\rvert}

\renewcommand{\d}{\textrm{d}}

\newcommand{\bbRp}{\bbR_{\geq 0}}
\newcommand{\bbRpp}{\bbR_{> 0}}

\newcommand{\bbZp}{\bbZ_{\geq 0}}
\newcommand{\bbZpp}{\bbZ_{> 0}}

\newcommand{\U}{\calU}
\newcommand{\Y}{\calY}

\newcommand{\Up}{\calU^-}
\newcommand{\Yf}{\calY^+}
\newcommand{\up}{u^-}
\newcommand{\yf}{y^+}

\newtheorem{definition}{Definition}
\newtheorem{lemma}{Lemma}
\newtheorem{theorem}{Theorem}

\newtheorem{remark}{Remark}

\newenvironment{myproof}{\proof}{\endproof}

\usetikzlibrary{arrows,automata, patterns, calc,decorations.pathmorphing,decorations.markings}

\tikzset{rbox/.style = {shape = rectangle, rounded corners,
		color=black,
		fill=yellow!20!white,
		text = black,
		inner sep = 1mm,
		outer sep= 0pt,
		minimum width = 12mm,
		minimum height = 6mm,
		draw}
}
\tikzset{box/.style = {shape = rectangle,
		color=black,
		fill=white!94!black,
		text = black,
		inner sep = 1mm,
		minimum width = 12mm,
		minimum height = 6mm,
		draw}
}
\tikzset{dbox/.style = {shape = rectangle,
		draw,
		dashed,
		fill=white!98!black}
}
\tikzset{drbox/.style = {shape = rectangle,
		rounded corners,
		draw,
		dashed,
		fill=yellow!10!white}
}
\tikzset{->,>=stealth',
	auto,
	node distance=3cm,
	semithick
}

\title{\LARGE \bfseries
	Completely monotone memory systems$^\ast$
}

\author{Brayan M. Shali$^1$, \quad Rodolphe Sepulchre$^{1,2}$
	\thanks{$^\ast$The research leading to these results has received funding from the European Research Council under the Advanced ERC Grant Agreement SpikyControl n.101054323. Email: {\itshape brayan.shali@kuleuven.be}, {\itshape rodolphe.sepulchre@kuleuven.be}.}%
	\thanks{$^1$Department of Electrical Engineering (ESAT), KU Leuven, KasteelPark Arenberg 10, B-3001 Leuven, Belgium.
	}
	\thanks{$^2$Department of Engineering, University of Cambridge, Trumpington Street, Cambridge CB2 1PZ, United Kingdom.}%
}

\begin{document}
	\maketitle
	\thispagestyle{empty}
	\pagestyle{empty}
	
	\begin{abstract}
		Linear time-invariant (LTI) relaxation systems are characterized by a completely monotone impulse response. This paper proposes an extension of this property to nonlinear systems: a system is said to have completely monotone memory if a positive increment in the past input leads to a completely monotone increment in the future output. We show that this property constitutes a strong form of monotonicity, while also imposing stability constraints that closely resemble fading memory, contraction, and incremental stability. We also show that completely monotone memory is equivalent to classical relaxation for LTI systems and to a combination of monotonicity and stability for a class of scalar state-affine systems. Finally, we illustrate the relevance of our framework to the gating dynamics of conductance-based models.
	\end{abstract}
	
	\section{Introduction}
	The concept of relaxation has its origins in the study of physical systems with memory, most notably in the theory of linear viscoelasticity \cite{gurtin1962}, where the stress at a given time depends on the history of the deformation. In this setting, the relaxation modulus describes how stress decays over time following a sudden deformation that is subsequently held fixed. In an important class of linear viscoelastic materials, the relaxation modulus admits a representation as a positive superposition of exponential decays \cite{husain2005, hanyga2013}. By Bernstein's theorem \cite{bernstein1929, widder1941}, this representation is equivalent to complete monotonicity of the relaxation modulus. This property characterizes a particular form of relaxation in which the response decays without oscillations, which is also the form of relaxation traditionally studied in systems theory.
	
	For linear time-invariant (LTI) systems, relaxation is classically characterized by complete monotonicity of the impulse response \cite{willems1972b}. This property admits several equivalent characterizations, including cyclic monotonicity of the Hankel operator \cite{chaffey2023}, variation diminution of the Hankel operator \cite{grussler2017}, gradient state-space realizations \cite{vanderschaft2024}, and RC/RL circuit realizations \cite{willems1976}. A broad perspective on relaxation systems and their connections to network and systems theory is provided in \cite{hancock2025}. The particular structure of relaxation systems has also been exploited in optimal control \cite{pates2019} and data-driven control \cite{drummond2024}. More recently, the theory has been extended to linear time-and-space-invariant (LTSI) systems \cite{donchev2025}, where it is shown that its core properties extend beyond the finite-dimensional LTI setting.
	
	The definition of relaxation in terms of a completely monotone impulse response cannot be directly extended to nonlinear systems. Existing attempts to generalize relaxation have therefore relied on alternative characterizations from the LTI setting. An input-output approach builds on cyclic monotonicity by requiring the memory functional of the system to be derived from a convex potential \cite{sepulchre2024}. However, this condition does not directly establish cyclic monotonicity of the nonlinear Hankel operator, which limits its interpretation as a direct extension of classical relaxation. A different, state-space approach extends the gradient realization perspective by defining nonlinear relaxation systems as passive Riemannian gradient systems \cite{vanderschaft2024}. In this formulation, relaxation is defined through the structure of a state-space realization, rather than directly in terms of input-output behavior.
	
	In this paper, we explore a different angle of attack that retains complete monotonicity as the central concept. Rather than extending one of the alternative characterizations of LTI relaxation, we formulate the property directly in terms of the Hankel operator, which maps past inputs to future outputs. Specifically, we introduce the notion of \emph{completely monotone memory}, which requires positive increments in past inputs to produce completely monotone increments in future outputs. This provides an input-output formulation of nonlinear relaxation that reduces to the classical definition in the LTI setting. We show that completely monotone memory constitutes a strong form of monotonicity \cite{angeli2003}, while also imposing stability constraints closely related to fading memory \cite{boyd1985}, contraction \cite{lohmiller1998}, and incremental stability \cite{angeli2002}. Finally, we illustrate the relevance of the framework to the gating dynamics of conductance-based models \cite{ermentrout2010}, while using the Hodgkin-Huxley model \cite{hodgkin1952} as a concrete example. In doing so, we establish that completely monotone memory is equivalent to a combination of monotonicity and stability for a class of scalar state-affine systems.
	
	The remainder of the paper is organized as follows. Section~\ref{sec:example} introduces the Hodgkin-Huxley model as an illustrative example. Section~\ref{sec:np} introduces notation and recalls the necessary preliminaries. Section~\ref{sec:ir} introduces the notion of completely monotone memory, establishes its relation to classical LTI relaxation, and discusses its connections to other system properties. Section~\ref{sec:cb_models} investigates completely monotone memory in the gating dynamics of conductance-based models. Section~\ref{sec:conclusion} concludes the paper.
	
	\section{Motivating example}\label{sec:example}
	
\begin{figure}
		\centering
		\begin{circuitikz}[-, /tikz/circuitikz/bipoles/length=30pt]
			\tikzset{every node/.style={font=\normalsize}}
			
			\draw[o->] (1,2) --node[midway,above] {$I_{\mathrm{ext}}$} (2,2);
			\draw[-] (1.5,2) -- (8,2);
			\draw[-] (8,0) -- (1.5,0);
			\draw[>-o] (2,0) -- (1,0);

			\draw (2.75,2) to[C, label=$C$] (2.75,0);
			
			\draw (4.5,2)
			to[vR, label=$g_{\mathrm{Na}}$] (4.5,.5)
			to[battery1, label=$E_{\mathrm{Na}}$] (4.5,0);
			
			\draw (6.25,2)
			to[vR, label=$g_{\mathrm{K}}$] (6.25,.5)
			to[battery1, label=$E_{\mathrm{K}}$] (6.25,0);
			
			\draw (8,2)
			to[R, label=$g_{\mathrm{L}}$] (8,.5)
			to[battery1, label=$E_{\mathrm{L}}$] (8,0);

		\end{circuitikz}
		\caption{Circuit representation of the Hodgkin-Huxley model.}
		\label{fig:HH}
	\end{figure}
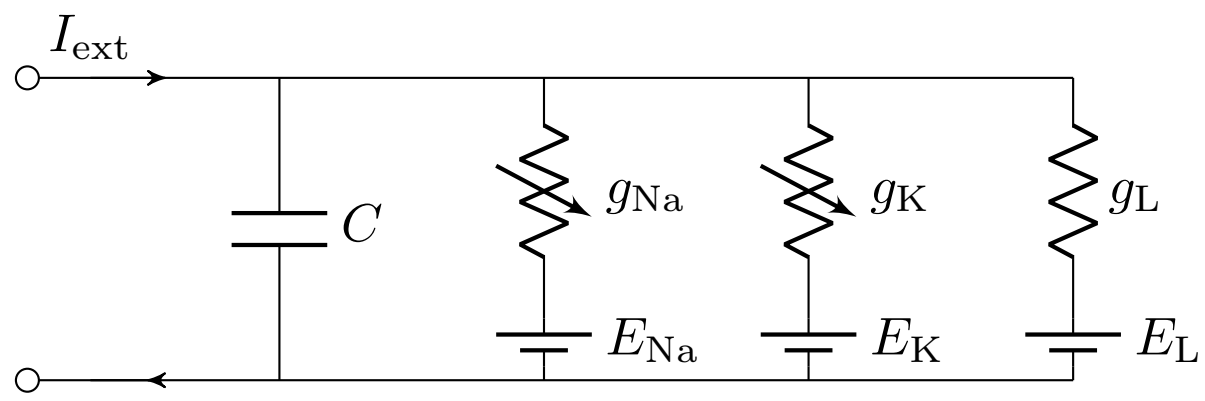

	The Hodgkin-Huxley model \cite{hodgkin1952}, see also \cite{ermentrout2010}, describes how the voltage across a neuron's membrane changes over time. Its equivalent electrical circuit is shown in Figure~\ref{fig:HH}, where a capacitor representing the membrane capacitance is connected in parallel with three branches representing the sodium, potassium, and leak currents. The membrane voltage dynamics follow the current balance equation
	
	\begin{equation}
		C \frac{dV}{dt}
		=
		I_{\mathrm{ext}}
		-I_{\mathrm{Na}}
		-I_{\mathrm{K}}
		-I_{\mathrm{L}}.
	\end{equation}
	Here, $V$ is the membrane voltage, $C$ is the membrane capacitance, and $I_{\mathrm{ext}}$ is the external current. The quantities $I_{\mathrm{Na}}$, $I_{\mathrm{K}}$, and $I_{\mathrm{L}}$ are the sodium, potassium, and leak currents, respectively. The currents are given by
	\begin{align}
		I_{\mathrm{Na}}	&= \bar g_{\mathrm{Na}}m^3h(V-E_{\mathrm{Na}}),\\
		I_{\mathrm{K}}	&= \bar g_{\mathrm{K}}n^4(V-E_{\mathrm{K}}),\\
		I_{\mathrm{L}}	&= g_{\mathrm{L}}(V-E_{\mathrm{L}}).
	\end{align}
	The quantities $\bar g_{\mathrm{Na}}$ and $\bar g_{\mathrm{K}}$ are the maximum sodium and potassium conductances, and $g_{\mathrm{L}}$ is the leak conductance. The quantities $E_{\mathrm{Na}}$, $E_{\mathrm{K}}$, and $E_{\mathrm{L}}$ are the reversal potentials. The variables $m$, $h$, and $n$ are gating variables that describe the opening and closing of voltage-dependent ion channels. 
	
	The dynamics of the gating variables are given by	
	\begin{align}
		\dot m &= \alpha_m(V)(1-m)-\beta_m(V)m,\\
		\dot h &= \alpha_h(V)(1-h)-\beta_h(V)h,\\
		\dot n &= \alpha_n(V)(1-n)-\beta_n(V)n,
	\end{align}
	where the functions $\alpha_{\ast}(V)$ and $\beta_{\ast}(V)$ are
	voltage-dependent transition rates given by
	\begin{align*}
		\alpha_m(V) &= \dfrac{0.1(V+40)}{1-e^{-(V+40)/10}},	& \beta_m(V)  &= 4e^{-(V+65)/18},\\
		\alpha_h(V) &= 0.07e^{-(V+65)/20}, 					& \beta_h(V)  &= \dfrac{1}{1+e^{-(V+35)/10}},\\
		\alpha_n(V) &= \dfrac{0.01(V+55)}{1-e^{-(V+55)/10}},& \beta_n(V)  &= 0.125e^{-(V+65)/80}.
	\end{align*}
	The gating dynamics are monotone \cite{angeli2003} with respect to the voltage input. More precisely, for the activation variables $m$ and $n$, an ordering of the voltage inputs and initial gating states is preserved by the corresponding trajectories. For the inactivation variable $h$, the ordering with respect to the voltage input is reversed, so that the dynamics of $h$ are more accurately described as ``antimonotone''. These properties follow from the signs of the voltage derivatives of the corresponding vector fields, see \cite[Proposition~III.2]{angeli2003} for details. In fact, the gating dynamics exhibit a stronger relaxation property.
	
	Specifically, consider the gating variable $m$. Consider two voltage inputs \(V_1\) and \(V_2\) and their corresponding gating trajectories $m_1$ and $m_2$. Suppose that $m_1(-T)= m_2(-T)$ and $V_1(t)\geq V_2(t)$ for $t\in[-T,0]$. By monotonicity of the dynamics with respect to the voltage input, it follows that
	\begin{equation}\label{eq:m10>m20}
		m_1(0)\ge m_2(0).
	\end{equation}
	Now, suppose that $V_1(t)=V_2(t)=0$ for all $t\ge 0$. Since both $m_1$ and $m_2$ evolve under the same constant voltage, subtracting their governing equations yields
	\begin{equation}
		\dot m_1(t) - \dot m_2(t)
		=-\lambda_m \left(m_1(t) - m_2(t)\right)
	\end{equation}
	for all $t\geq 0$, where
	\begin{equation}
		\lambda_m=\alpha_m(0)+\beta_m(0)>0.
	\end{equation}
	Consequently, we have that
	\begin{equation}
		m_1(t) - m_2(t)=\left(m_1(0) - m_2(0)\right)e^{-\lambda_m t},
	\end{equation}
	which, due to \eqref{eq:m10>m20}, is a nonnegative decaying exponential. 
	
	This can be viewed as an incremental version of the classical relaxation property for LTI systems. Whereas relaxation requires the future response to a nonnegative past perturbation to be completely monotone, the incremental version concerns the difference between future responses corresponding to ordered past inputs. In the present example, this difference is nonnegative and decays to zero in a completely monotone manner once both voltage inputs are set to zero. Thus, the difference between the two trajectories relaxes without oscillations.
	
	In the following sections, we formalize this observation through the concept of completely monotone memory. We then investigate its relation to classical relaxation and monotonicity, and return to the Hodgkin-Huxley gating dynamics to establish the property rigorously.
	
	\section{Notation and preliminaries}\label{sec:np}
	We denote the set of real numbers (integers) by $\bbR$ ($\bbZ$). We denote the sets of nonnegative and positive real numbers (integers) by $\bbRp$ ($\bbZp$) and $\bbRpp$ ($\bbZpp$), respectively. We use similar notation for the sets of nonpositive and negative real numbers (integers). We say that a function $f:X\to\bbR$ is nonnegative if $f(x) \geq 0$ for all $x\in X$. Positive, nonpositive, and negative functions are defined in a similar manner.

	\subsection{Completely monotone functions}\label{sec:np_cm}
	
	A function $f:\bbRp\to\bbR$ is \emph{completely monotone} if it is smooth and
	\begin{equation}
		(-1)^k f^{(k)}(x)\geq0,\qquad k\in\bbZp,\quad x\in\bbRp.
	\end{equation}
	By Bernstein's theorem \cite{bernstein1929}, see also \cite{widder1941, merkle2014}, a function is completely monotone if and only if it admits a representation
	\begin{equation}
		f(x)=\int_0^\infty e^{-\lambda x}\,\nu(\d\lambda),
		\qquad x\in\bbRp,
	\end{equation}
	where $\nu$ is a positive Borel measure on $\bbRp$ with finite moments of every order. Thus, completely monotone functions can be viewed as positive mixtures of decaying exponentials. The sum and product of two completely monotone functions are also completely monotone. In particular, this means that completely monotone functions form a convex cone.

	\subsection{Systems as operators}	
	We view systems as operators $G:\U\to\Y$ that map input signals to output signals. The choice of signal spaces $\U$ and $\Y$ affects the class of systems under consideration. For example, requiring bounded outputs for bounded inputs excludes the ideal integrator, whereas requiring bounded outputs for absolutely integrable inputs includes it. There is no unique natural choice, and we do not pursue this issue further. For concreteness, we take $\U$ to be the space of bounded piecewise-continuous signals and $\Y$ the space of bounded continuous signals, assuming throughout that the operators and properties considered are well-defined.
	
	With this in mind, for $t\in\bbR$, we define the time-shift operator $S_t:\U\to\U$ by
	\begin{equation}
		(S_tu)(\tau)=u(\tau+t),
		\qquad \tau\in\bbR.
	\end{equation}
	The operator $G$ is \emph{time-invariant} (TI) if
	\begin{equation}\label{eq:G_TI}
		G(S_tu)=S_tG(u)
	\end{equation}
	for all $t\in\bbR$ and $u\in\U$. In other words, shifting the input in time results in the same shift of the output. The operator $G$ is \emph{causal} if, for all $t\in\bbR$ and $u,v\in\U$, $u(\tau)=v(\tau)$ for all $\tau\leq t$ implies that $G(u)(t)=G(v)(t)$. Thus, the output at any time depends only on the input up to that time.
	
	We associate a causal TI operator $G$ with a \emph{memory functional} $M$, which maps past input signals to the present output. To define this functional, let $P:\U\to\Up$ denote the past truncation operator given by
	\begin{equation}
		(Pu)(t)=u(t),\qquad t\leq 0,
	\end{equation}
	where $\Up=P\U$ is the space of past input signals obtained by truncating signals in $\U$ at time zero. Let $E:\Up\to\U$ denote the zero-extension operator given by
	\begin{equation}
		(E\up)(t)=
		\begin{dcases}
			\up(t), & t\leq 0,\\
			0, & t>0.
		\end{dcases}
	\end{equation}
	The memory functional $M:\Up\to\bbR$ is then defined by
	\begin{equation}
		M(\up)=G(E\up)(0).
	\end{equation}
	Causality and time invariance imply that
	\begin{equation}
		G(u)(t)=M(PS_tu),\qquad t\in\bbR,
	\end{equation}
	hence the memory functional $M$ uniquely determines $G$.
	
	We also associate $G$ with a \emph{Hankel operator} $H$, which maps past input signals to future output signals. To define this operator, let $F:\Y\to\Yf$ denote the future truncation operator given by
	\begin{equation}
		(Fy)(t)=y(t),\qquad t\geq0,
	\end{equation}
	where $\Yf=F\Y$ is the space of future output signals. The Hankel operator $H:\Up\to\Yf$ is defined by
	\begin{equation}\label{eq:H_definition}
		H(\up)=FG(E\up).
	\end{equation}
	By construction, we have that
	\begin{equation}
		M(\up)=H(\up)(0),
	\end{equation}
	hence $H$ uniquely determines $G$ through its memory functional $M$. We note that not every $H:\Up\to\Yf$ is the Hankel operator of a causal TI operator $G:\U\to\Y$. In particular, $H$ must satisfy the consistency condition
	\begin{equation}\label{eq:hankel_consistency}
		H(PS_tE\up)(\tau)=H(\up)(t+\tau),
		\qquad t,\tau\geq0.
	\end{equation}
	Indeed, $(S_tE\up)(\tau) = 0$ for $t\geq0$ and $\tau > 0$, so
	\begin{equation}
		EPS_tE\up=S_tE\up,\qquad t\geq 0.
	\end{equation}
	The condition then follows from  \eqref{eq:H_definition} and \eqref{eq:G_TI}.
	
	\section{Completely monotone memory}\label{sec:ir}
	A defining characteristic of relaxation is that a nonnegative perturbation in the past produces a completely monotone response in the future. For causal LTI systems, this property means complete monotonicity of the impulse response \cite{willems1972b}. Since the concept of impulse response is problematic for nonlinear systems, we instead formulate relaxation explicitly as a property of the system's memory, expressed via the associated Hankel operator. We also take an incremental perspective on relaxation, leading to the following definition.
	
	\begin{definition}
		A causal TI operator $G:\U\to\Y$ has \emph{completely monotone memory} if the associated Hankel operator $H:\Up\to\Yf$ is such that $H(\up_1)-H(\up_2)$ is completely monotone whenever $\up_1-\up_2$ is nonnegative.
	\end{definition}
	
	The remainder of this section is split into two parts. First, we establish the equivalence between completely monotone memory and relaxation for LTI systems. Second, we investigate the relation between completely monotone memory and monotonicity, and discuss its connections to fading memory, contraction, and incremental stability.
	
	\subsection{Linear systems with completely monotone memory}
	
	We first consider the special case of linear systems, for which the definition of completely monotone memory admits a simpler characterization. In particular, linearity allows the condition on differences of Hankel outputs to be expressed in terms of the response to a single nonnegative past input.
	\begin{lemma}\label{lem:LTI_equivalence}
		Let $G:\U\to\Y$ be a causal LTI operator with associated Hankel operator $H:\Up\to\Yf$.	Then, $G$ has completely monotone memory if and only if $H(\up)$ is completely monotone whenever $\up$ is nonnegative.
	\end{lemma}
	\begin{myproof}
		Since $G$ is linear, its Hankel operator $H$ is also linear, so that
		\begin{equation}
			H(\up_1)-H(\up_2)=H(\up_1-\up_2)
		\end{equation}
		for all $\up_1,\up_2\in\Up$. Therefore, the condition for completely monotone memory is equivalent to requiring that $H(\up)$ be completely monotone whenever $\up$ is nonnegative.
	\end{myproof}
	\begin{remark}
		The condition that $H(\up)$ is completely monotone whenever $\up$ is nonnegative can be interpreted as a non-incremental form of completely monotone memory. It requires the future response to a nonnegative past input to be completely monotone, whereas completely monotone memory imposes this on differences between future responses to ordered past inputs. Lemma~\ref{lem:LTI_equivalence} establishes the equivalence between both forms for linear systems. This is not the case in general. In particular, the incremental form implies the non-incremental form whenever \(H(0)\) is completely monotone, but the converse is not true.
	\end{remark}
	
	For LTI systems, Willems \cite{willems1972b,willems1976} formalized the notion of relaxation as complete monotonicity of the impulse response. He established connections between relaxation and properties such as passivity, reciprocity, and the existence of so-called compatible storage functions. More recently, relaxation has been characterized in terms of cyclic monotonicity of the associated Hankel operator \cite{chaffey2023}. Although cyclic monotonicity is distinct from the property imposed by our definition of completely monotone memory, we show that both characterize the same class of LTI convolution systems.
	
	To this end, consider a causal LTI system with input-output operator \(G:\U\to\Y\) given by
	\begin{equation}\label{eq:G_lti}
		G(u)(t)= \int_{-\infty}^t g(t-\tau)u(\tau)\,\d\tau,
	\end{equation}
	where $g:\bbRp\to\bbR$ is the impulse response. The associated Hankel operator $H:\Up\to\Yf$ is given by
	\begin{equation}\label{eq:H_lti}
		H(\up)(t)
		=\int_{-\infty}^0 g(t-\tau)\up(\tau)\, \d\tau,
		\qquad t\geq0.
	\end{equation}
	We implicitly assume that $g$ is such that $G$ and $H$ are well-defined. With $\U$ the space of bounded piecewise continuous signals and $\Y$ the space of bounded continuous signals, it is sufficient to assume that $g$ is absolutely integrable.
	
	The following theorem shows the equivalence between completely monotone memory and relaxation for LTI systems. As such, it supports the interpretation of completely monotone memory as a nonlinear extension of relaxation.	
	\begin{theorem}\label{thm:LTI_equivalence}
		Consider an LTI system with impulse response $g:\bbRp\to\bbR$ and input-output operator $G:\calU\to\calY$ given by \eqref{eq:G_lti}. Suppose that $g$ is continuously differentiable, and that both $g$ and its derivative are absolutely integrable. Then, $G$ has completely monotone memory if and only if $g$ is completely monotone.
	\end{theorem}
	\begin{myproof}
			First, suppose that $g$ is completely monotone. By Bernstein's	theorem \cite{bernstein1929, widder1941}, we can write
			\begin{equation}\label{eq:g_Bernstein}
				g(t)=\int_0^\infty e^{-\alpha t}\,\nu(\d\alpha),
			\end{equation}
			where $\nu$ is a positive Borel measure on $\bbRp$ with finite moments of every order. The associated Hankel operator $H:\Up\to\Yf$ is given by \eqref{eq:H_lti}, where the change of variables $s=-\tau$ yields
			\begin{equation}\label{eq:H_lti_flipped}
				H(\up)(t)=\int_0^\infty g(t+s)\up(-s)\,\d s.
			\end{equation}
			Let $\up\in\Up$ be a nonnegative past input. Substituting
			\eqref{eq:g_Bernstein} in \eqref{eq:H_lti_flipped} and applying
			Tonelli's theorem yields
			\begin{align*}
				H(\up)(t)
				&=\int_0^\infty\int_0^\infty
				e^{-\alpha t}e^{-\alpha s}\up(-s)
				\,\d s\,\nu(\d\alpha)\\
				&=\int_0^\infty e^{-\alpha t}
				\left(\int_0^\infty e^{-\alpha s}\up(-s)\,\d s\right)
				\nu(\d\alpha).
			\end{align*}
			The measure
			\begin{equation}
				\mu(\d\alpha)=\left(\int_0^\infty e^{-\alpha s}\up(-s)\,\d s\right)\nu(\d\alpha)
			\end{equation}
			is positive with finite moments of every order, hence $H(\up)$ is completely monotone by Bernstein's theorem. Therefore, by Lemma~\ref{lem:LTI_equivalence}, $G$ has completely monotone memory.
			
			Conversely, suppose that $G$ has completely monotone memory. For $\lambda>0$, let the nonnegative past input
			$\up_\lambda\in\Up$ be given by $\up_\lambda(t)=e^{\lambda t}$, and define $H_\lambda=H(\up_\lambda)$, so that
			\begin{equation}
				H_\lambda(t)=\int_0^\infty
				g(t+s)e^{-\lambda s}\,\d s.
			\end{equation}
			By Lemma~\ref{lem:LTI_equivalence}, $H_\lambda$ is completely
			monotone for every $\lambda>0$. Since $g$ and its
			derivative are absolutely integrable, we can
			differentiate under the integral to obtain
			\begin{equation}
				\dot H_\lambda(t)
				=\int_0^\infty g'(t+s)e^{-\lambda s}\,\d s
				=-g(t)+\lambda H_\lambda(t).
			\end{equation}
			Both $\lambda H_\lambda$ and $-\dot H_\lambda$ are
			completely monotone, hence
			\[
			g=\lambda H_\lambda-\dot H_\lambda
			\]
			is also completely monotone.
	\end{myproof}
	
	\subsection{Completely monotone memory and monotonicity}
	
	We next investigate the relation between completely monotone memory and monotonicity, as well as its connections to fading memory, contraction, and incremental stability. To begin with, completely monotone memory can be viewed as an order preservation property. Indeed, the convex cone of nonnegative functions induces a partial order $\geq$ on $\Up$, defined by
	\begin{equation}
		u_1 \geq u_2
		\quad\Longleftrightarrow\quad
		u_1-u_2 \text{ is nonnegative}.
	\end{equation}
	Similarly, the convex cone of completely monotone functions induces a partial order $\geq_{\mathrm{cm}}$ on $\Yf$, defined by
	\begin{equation}
		y_1 \geq_{\mathrm{cm}} y_2
		\quad\Longleftrightarrow\quad
		y_1-y_2 \text{ is completely monotone}.
	\end{equation}
	Thus, completely monotone memory can be formulated as
	\begin{equation}
		\up_1\geq \up_2 \quad\implies\quad H(\up_1)\geq_{\mathrm{cm}}H(\up_2),
	\end{equation}
	which states that the Hankel operator is order preserving with respect to the orders above. 
	
	An important consequence of this interpretation is that completely monotone memory implies monotonicity \cite{angeli2003}. Indeed, since every completely monotone function is nonnegative, completely monotone memory implies that
	\begin{equation}
		\up_1\geq \up_2 \quad\implies\quad H(\up_1)\geq H(\up_2).
	\end{equation}
	Evaluating the Hankel output at time $0$ yields
	\begin{equation}
		\up_1\geq \up_2 \quad\implies\quad M(\up_1)\geq M(\up_2),
	\end{equation}
	where $M$ is the associated memory functional. Therefore, by causality and time-invariance, if $G$ has completely monotone memory, then it is also monotone, in the sense that
	\begin{equation}
		u_1\ge u_2\quad\implies\quad G(u_1)\ge G(u_2).
	\end{equation}
	Therefore, completely monotone memory strengthens monotonicity by requiring the difference between the future outputs generated by ordered past inputs to be completely monotone, rather than merely nonnegative. In particular, this difference is nonincreasing, so that the influence of a nonnegative perturbation in the past cannot grow once the input is set to zero. If the difference additionally converges to zero, the effect of the perturbation vanishes asymptotically. Although this convergence is not guaranteed by complete monotonicity alone, it may be imposed as an additional requirement, leading to what we may call \emph{vanishing completely monotone memory}. This amounts to requiring that the output differences be completely monotone functions vanishing at infinity, or equivalently, by Bernstein's theorem, that their representing measures have no atom at zero.
	
	Vanishing completely monotone memory is reminiscent of fading memory \cite{boyd1985}. To make the comparison more explicit, we focus on a particular consequence of fading memory: for any two past inputs and any common future input, the corresponding future outputs converge to each other as time tends to infinity. Intuitively, this follows because the differences between the past inputs become increasingly remote and therefore have a diminishing influence on the outputs. Vanishing completely monotone memory shares this asymptotic convergence property, but only for \emph{ordered} past inputs under a common \emph{zero} future input. In return, it imposes the stronger requirement that the output differences be completely monotone. Moreover, fading memory is fundamentally a continuity property with respect to a weighted norm on input histories, whereas vanishing completely monotone memory imposes no such continuity requirement. This raises the question of whether vanishing completely monotone memory, together with a weaker form of continuity, is sufficient to imply fading memory. The answer is not straightforward, and we leave a more detailed investigation for future work.
	
	The additional restrictions imposed by (vanishing) completely monotone memory can also be examined at the level of the system's dynamics. To illustrate this, consider a scalar state-space system of the form
	\begin{equation}
		\dot y=f(y,u),
	\end{equation}
	where \(f:\bbR\times\bbR\to\bbR\) is continuously differentiable and the dynamics are assumed to induce a well-defined causal TI input-output operator \(G:\U\to\Y\). For such a system, monotonicity is characterized by the condition 
	\begin{equation}\label{eq:f_u}
		\frac{\partial f}{\partial u}(y,u)\geq0,
	\end{equation}
	see \cite[Proposition~III.2]{angeli2003}. This ensures that ordered inputs generate ordered outputs, but completely monotone memory imposes additional restrictions on the free response.
	
	To see this, consider two ordered past inputs \(\up_1\geq \up_2\). Let $\yf_1=H(\up_1)$ and $\yf_2=H(\up_2)$ denote the corresponding future outputs. By completely monotone memory, $\yf_1-\yf_2$ is completely monotone, hence, in particular,
	\begin{equation}
		\yf_1(t) - \yf_2(t) \geq 0,
		\qquad
		\dot{y}^+_1(t)-\dot{y}^+_2(t) \leq 0.
	\end{equation}
	Since the future input is zero, the future outputs satisfy
	\begin{equation}
		\dot{y}^+_1=f(\yf_1,0),\qquad \dot{y}^+_2=f(\yf_2,0),
	\end{equation}
	so that
	\begin{equation}
		f(\yf_1(t),0)-f(\yf_2(t),0)\leq0.
	\end{equation}
	Thus, whenever \(\yf_1(t) - \yf_2(t) \neq 0\), we have
	\begin{equation}
		\frac{f(\yf_1(t),0)-f(\yf_2(t),0)}
		{\yf_1(t)-\yf_2(t)}\leq0.
	\end{equation}
	Under suitable reachability assumptions, which we do not pursue here, this suggests the necessary condition
	\begin{equation}\label{eq:f_y}
			\frac{\partial f}{\partial y}(y,0)\leq0
	\end{equation}
	for all reachable states $y$.
	
	Condition \eqref{eq:f_y} implies that the free response is nonexpansive on the relevant state domain, in the sense that the distance between any two trajectories cannot increase. In contraction theory, this property is referred to as \emph{weak contraction} \cite{bullo2026}, whereas \emph{contraction} \cite{lohmiller1998} follows from the stronger condition 
	\begin{equation}
		\frac{\partial f}{\partial y}(y,0)\leq-\alpha,\qquad \alpha > 0,
	\end{equation}
	which guarantees that trajectories converge to each other at an exponential rate. Another closely related notion is incremental asymptotic stability \cite{angeli2002}, which requires trajectories driven by the same input to converge to each other, without necessarily imposing an exponential rate of convergence.
	
	As in the comparison with fading memory, (vanishing) completely monotone memory differs from these notions in two important aspects. First, it only constrains the evolution of output differences generated by ordered past inputs when the future input is zero, rather than arbitrary trajectories driven by a common input. Second, it requires these differences to be completely monotone, which imposes restrictions on all their derivatives, rather than merely requiring nonexpansiveness or asymptotic convergence. In particular, completely monotone memory does not necessarily imply incremental asymptotic stability, while neither contraction nor incremental asymptotic stability is generally sufficient to guarantee completely monotone memory, even when combined with monotonicity.
	
	\begin{remark}\label{rem:unordered_inputs}
		Completely monotone memory bounds future output differences even when they are generated by \emph{unordered} past inputs. Indeed, given any $\up_1,\up_2\in\Up$, define their pointwise minimum $\ell=\min\{\up_1,\up_2\}$ and maximum $r=\max\{\up_1,\up_2\}$. Then $\ell\leq\up_i\leq r$ for $i\in\{1,2\}$, so monotonicity yields $H(\ell)\leq H(\up_i)\leq H(r)$ and, thus,
		\begin{equation}
			\abs{H(\up_1)(t)-H(\up_2)(t)} \leq H(r)(t)-H(\ell)(t).
		\end{equation}
		The upper bound is completely monotone and hence nonincreasing. For vanishing completely monotone memory, it converges to zero, so the output difference vanishes even for unordered past inputs, although not necessarily in completely monotone fashion.
	\end{remark}
	
	\section{Completely monotone memory in conductance-based models}\label{sec:cb_models}
	
	We finally investigate the property of completely monotone memory in the gating dynamics of conductance-based models. As introduced in Section~\ref{sec:example}, the three gating variables $m$, $h$, and $n$ of the Hodgkin-Huxley model share the same general form of dynamics. In particular, they are described by scalar state-affine systems of the form
	\begin{equation}\label{eq:ltv}
		\dot y=-a(u)y+b(u),
	\end{equation}
	where \(a:\bbR\to\bbR\) and \(b:\bbR\to\bbR\) are continuously differentiable. Throughout this section, we assume that $a$ is bounded below by a positive constant, as is the case for all three Hodgkin-Huxley gating variables. We emphasize that this structure is not specific to the Hodgkin--Huxley model, but that scalar differential equations of this form are the canonical way of modeling the dynamics of voltage-dependent activation and inactivation gates in conductance-based models \cite{ermentrout2010}. They can, in fact, be regarded as mean-field models of the stochastic reactions that gate the conductance at the molecular level \cite{destexhe1994}. These gates typically respond monotonically to changes in voltage and relax toward an equilibrium when the voltage is held constant. As discussed in the preceding section, monotonicity and stability are closely related to the notion of completely monotone memory. We now establish this connection formally for systems of the form \eqref{eq:ltv}.
	
	For a given input signal $u\in\U$, \eqref{eq:ltv} is a linear time-varying differential equation. Given an initial condition $y(t_0) = y_0$ at time $t_0\in\bbR$, the solution \cite{brockett2015} for $t\geq t_0$ is given by
	\begin{equation}\label{eq:ltv_solution}
		y(t)
		= \phi(t,t_0,u)y_0
		+ \int_{t_0}^{t}\phi(t,\tau,u)b(u(\tau))\,\d\tau,
	\end{equation}
	where
	\begin{equation}\label{eq:phi}
		\phi(t,\tau,u) = \exp\left(-\int_{\tau}^{t}a(u(s))\, \d s\right).
	\end{equation}
	The positive lower bound on $a$ ensures that 
	\begin{equation}
		\lim_{\tau\to-\infty} \phi(t,\tau,u) = 0,
	\end{equation}
	hence the contribution of the initial condition in \eqref{eq:ltv_solution} vanishes as $t_0\to-\infty$. In particular, taking the limit $t_0\to-\infty$ yields
	\begin{equation}\label{eq:ltv_operator}
		y(t)=\int_{-\infty}^{t}
		\phi(t,\tau,u)b(u(\tau))\,d\tau.
	\end{equation}
	Since $b$ is continuous and $u$ is piecewise continuous and bounded, the composition $b\circ u$ is piecewise continuous and bounded. Hence, the positive lower bound on $a$ implies that the integral in \eqref{eq:ltv_operator} converges absolutely for all $t\in\bbR$, and the signal $y$ is continuous and bounded. Consequently, \eqref{eq:ltv} defines a causal TI operator $G:\U\to\Y$, where $y= G(u)$ is given by \eqref{eq:ltv_operator}. 
	
	To derive the associated Hankel operator $H:\Up\to\Yf$, we note that under zero input, the dynamics reduce to
	\begin{equation}
		\dot y=-a(0)y+b(0),
	\end{equation}
	so the future output is given by
	\begin{equation}
		y(t) = \frac{b(0)}{a(0)} + \left(y(0) - \frac{b(0)}{a(0)}\right)e^{-a(0)t}.
	\end{equation}
	In other words, the Hankel operator is given by
	\begin{equation}\label{eq:ltv_hankel_operator}
		H(\up)(t)
		=\frac{b(0)}{a(0)} + \left(M(\up) - \frac{b(0)}{a(0)}\right)e^{-a(0)t},
	\end{equation}
	where the memory functional $M:\Up\to\bbR$ is given by
	\begin{equation}\label{eq:memory_functional}
		M(\up) = \int_{-\infty}^{0}	\phi(0,\tau,\up)b(\up(\tau))\,\d\tau.
	\end{equation}
	This representation makes the completely monotone memory condition particularly simple. Given two past inputs $\up_1,\up_2\in\Up$, subtracting their Hankel outputs gives
	\begin{equation}\label{eq:ltv_Hankel_increment}
		H(\up_1)(t)-H(\up_2)(t)
		=\left(M(\up_1)-M(\up_2)\right)e^{-a(0)t}.
	\end{equation}
	Since $a(0)>0$, it follows that $H(\up_1)-H(\up_2)$ is completely monotone if and only if $M(\up_1)-M(\up_2)$ is nonnegative. Consequently, completely monotone memory can be established by showing that the system is monotone with respect to its input. We first illustrate this for the \(m\)-gate of the Hodgkin-Huxley model.
	
	Recall that the $m$-gate dynamics are given by
	\begin{equation}
		a(u)=\alpha_m(u)+\beta_m(u), \qquad b(u)=\alpha_m(u),
	\end{equation}
	where the rate functions $\alpha_m$ and $\beta_m$ are given in Section~\ref{sec:example}. The interval $[0,1]$ is forward invariant under these dynamics. Indeed, writing
	\begin{equation*}
		f(y,u) = -a(u)y + b(u) = \alpha_m(u)(1-y)-\beta_m(u)y
	\end{equation*}
	we have $f(0,u)=\alpha_m(u)>0$ and $f(1,u)=-\beta_m(u)<0$ for all $u\in\bbR$. Moreover, the dynamics are monotone with respect to the input on the interval $[0,1]$. Indeed, the rate functions satisfy $\alpha_m'(u)>0$ and $\beta_m'(u)<0$, hence
	\begin{equation*}
		\frac{\partial f}{\partial u}(y,u)= \alpha_m'(u)(1-y)-\beta_m'(u)y\geq0, \quad y\in[0,1].
	\end{equation*}		
	Thus, by \cite[Proposition~III.2]{angeli2003}, for any initial time $t_0$ and final time $T > t_0$, if $u_1(t)\geq u_2(t)$ for all $t\in[t_0,T]$, then the corresponding solutions $y_1$ and $y_2$, initialized at	the same state $y_1(t_0)=y_2(t_0)\in[0,1]$, satisfy
	\begin{equation}
		y_1(t)\geq y_2(t),\qquad t\in[t_0,T].
	\end{equation}
	Now, fixing $T = 0$ and taking the limit $t_0 \to-\infty$ shows that $M(\up_1)- M(\up_2)$ is nonnegative whenever $\up_1 - \up_2$ is nonnegative. In view of \eqref{eq:ltv_Hankel_increment}, we conclude that the $m$-gate dynamics have completely monotone memory.
	
	Since the rate functions of the $n$-gate satisfy $\alpha_n'(u)>0$ and $\beta_n'(u)<0$, the same arguments show that the {$n$-gate} dynamics also have completely monotone memory. In contrast, the rate functions of the $h$-gate satisfy $\alpha_h'(u)<0$ and $\beta_h'(u)>0$, hence the dynamics are monotone in the opposite direction with respect to the input. This means that the $h$-gate dynamics have reversed completely monotone memory, that is, $H(\up_1)-H(\up_2)$ is completely monotone whenever $\up_1-\up_2$ is \emph{nonpositive}.
	
	The preceding arguments rely only on the stability of the free response and the monotonicity of the system with respect to its input. Therefore, they can be summarized in the following theorem, which establishes the equivalence between monotonicity and completely monotone memory for this class of stable state-affine systems.

	\begin{theorem}
		Consider the system \eqref{eq:ltv}, where \(a:\bbR\to\bbR\) and \(b:\bbR\to\bbR\) are continuously differentiable. Suppose that $a$ is bounded from below by a positive constant. Then, \eqref{eq:ltv} defines a causal TI operator \(G:\U\to\Y\) given by
		\begin{equation}
			G(u)(t)=\int_{-\infty}^t
			\phi(t,\tau,u)b(u(\tau))\,\mathrm d\tau,
		\end{equation}
		where \(\phi(t,\tau,u)\) is given by \eqref{eq:phi}. Furthermore, $G$ has completely monotone memory if and only if it is monotone.
	\end{theorem}
	
	This theorem relies on the particular structure of the free response dynamics, which are scalar and linear. In this case, stability ensures that differences between future outputs decay exponentially, so that monotonicity is sufficient to guarantee completely monotone memory. This equivalence does not generally extend to higher-dimensional systems, nonlinear output maps, or general nonlinear dynamics, where stability and monotonicity alone are insufficient to guarantee complete monotonicity of the output differences.
	
	\section{Conclusion}\label{sec:conclusion}
	
	We introduced completely monotone memory as an extension of the classical relaxation property to nonlinear systems that retains complete monotonicity as the central concept. By formulating the property in terms of the Hankel operator, we showed that completely monotone memory constitutes a strong form of monotonicity and discussed its connections to fading memory, contraction, and incremental stability. We also established its equivalence to classical relaxation for LTI systems and to a combination of stability and monotonicity for a class of scalar state-affine systems. The latter result was illustrated through the gating dynamics of the Hodgkin-Huxley model. 
	
	It is of particular interest that the proposed strong form of monotonicity also imposes stability constraints on the free response. These two properties are usually imposed separately. It is also noteworthy that the gating dynamics of conductance-based models have this property. These dynamics model memory mechanisms in neuronal circuits, which suggests an important form of relaxation to be acknowledged in memristive modeling \cite{forni2025} and neuromorphic design.
	\bibliographystyle{ieeetr}
	\bibliography{../../references/kernel_refs}	

\begin{thebibliography}{10}

\bibitem{gurtin1962}
M.~E. Gurtin and E.~Sternberg, ``On the linear theory of viscoelasticity,''
  {\em Archive for Rational Mechanics and Analysis}, vol.~11, no.~1,
  pp.~291--356, 1962.

\bibitem{husain2005}
S.~A. Husain and R.~Anderssen, ``Modelling the relaxation modulus of linear
  viscoelasticity using kohlrausch functions,'' {\em Journal of Non-Newtonian
  Fluid Mechanics}, vol.~125, no.~2-3, pp.~159--170, 2005.

\bibitem{hanyga2013}
A.~Hanyga, ``Wave propagation in linear viscoelastic media with completely
  monotonic relaxation moduli,'' {\em Wave Motion}, vol.~50, no.~5,
  pp.~909--928, 2013.

\bibitem{bernstein1929}
S.~Bernstein, ``Sur les fonctions absolument monotones,'' {\em Acta
  Mathematica}, vol.~52, no.~0, pp.~1--66, 1929.

\bibitem{widder1941}
D.~V. Widder, {\em The Laplace Transform}.
\newblock Princeton University Press, 1941.

\bibitem{willems1972b}
J.~C. Willems, ``{Dissipative dynamical systems Part II: Linear systems with
  quadratic supply rates},'' {\em Archive for Rational Mechanics and Analysis},
  vol.~45, no.~5, pp.~352--393, 1972.

\bibitem{chaffey2023}
T.~Chaffey, H.~J. van Waarde, and R.~Sepulchre, ``Relaxation systems and cyclic
  monotonicity,'' in {\em Proceedings of the IEEE Conference on Decision and
  Control}, pp.~1673--1679, 2023.

\bibitem{grussler2017}
C.~Grussler, J.~Umenberger, and I.~R. Manchester, ``Identification of
  externally positive systems,'' in {\em Proceedings of the IEEE Conference on
  Decision and Control}, pp.~6549--6554, 2017.

\bibitem{vanderschaft2024}
A.~van~der Schaft, ``Reciprocity of nonlinear systems,'' {\em SIAM Journal on
  Control and Optimization}, vol.~62, no.~6, pp.~3019--3041, 2024.

\bibitem{willems1976}
J.~C. Willems, ``Realization of systems with internal passivity and symmetry
  constraints,'' {\em Journal of the Franklin Institute}, vol.~301, no.~6,
  pp.~605--621, 1976.

\bibitem{hancock2025}
A.~J. Hancock, ``On relaxation systems in network and systems theory,''
  {Master's Thesis}, University of Cambridge, 2025.

\bibitem{pates2019}
R.~Pates, C.~Bergeling, and A.~Rantzer, ``On the optimal control of relaxation
  systems,'' in {\em Proceedings of the IEEE Conference on Decision and
  Control}, pp.~6068--6073, 2019.

\bibitem{drummond2024}
R.~Drummond, H.~Taghavian, and P.~R. Baldivieso-Monasterios,
  ``Learning-to-control relaxation systems with the step response,'' in {\em
  Proceedings of the American Control Conference}, pp.~3584--3589, 2024.

\bibitem{donchev2025}
T.~I. Donchev, B.~M. Shali, and R.~Sepulchre, ``Linear time-and-space-invariant
  relaxation systems,'' in {\em Proceedings of the IEEE Conference on Decision
  and Control}, pp.~4113--4118, 2025.

\bibitem{sepulchre2024}
R.~Sepulchre, T.~Chaffey, F.~Forni, and Y.~Huo, ``On the relaxation property of
  nonlinear circuit elements,'' in {\em Proceedings of the IEEE Conference on
  Decision and Control}, 2024.

\bibitem{angeli2003}
D.~Angeli and E.~Sontag, ``Monotone control systems,'' {\em IEEE Transactions
  on Automatic Control}, vol.~48, pp.~1684--1698, Oct. 2003.

\bibitem{boyd1985}
S.~Boyd and L.~Chua, ``Fading memory and the problem of approximating nonlinear
  operators with volterra series,'' {\em IEEE Transactions on Circuits and
  Systems}, vol.~32, no.~11, pp.~1150--1161, 1985.

\bibitem{lohmiller1998}
W.~Lohmiller and J.-J.~E. Slotine, ``On contraction analysis for non-linear
  systems,'' {\em Automatica}, vol.~34, no.~6, pp.~683--696, 1998.

\bibitem{angeli2002}
D.~Angeli, ``A lyapunov approach to incremental stability properties,'' {\em
  IEEE Transactions on Automatic Control}, vol.~47, pp.~410--421, Mar. 2002.

\bibitem{ermentrout2010}
G.~B. Ermentrout and D.~H. Terman, {\em Mathematical Foundations of
  Neuroscience}.
\newblock Springer New York, 2010.

\bibitem{hodgkin1952}
A.~L. Hodgkin and A.~F. Huxley, ``A quantitative description of membrane
  current and its application to conduction and excitation in nerve,'' {\em The
  Journal of Physiology}, vol.~117, no.~4, pp.~500--544, 1952.

\bibitem{merkle2014}
M.~Merkle, {\em Completely Monotone Functions: A Digest}, pp.~347--364.
\newblock Springer New York, 2014.

\bibitem{bullo2026}
F.~Bullo, {\em Contraction Theory for Dynamical Systems}.
\newblock Kindle Direct Publishing, {1.3}~ed., 2026.

\bibitem{destexhe1994}
A.~Destexhe, Z.~F. Mainen, and T.~J. Sejnowski, ``Synthesis of models for
  excitable membranes, synaptic transmission and neuromodulation using a common
  kinetic formalism,'' {\em Journal of Computational Neuroscience}, vol.~1,
  no.~3, pp.~195--230, 1994.

\bibitem{brockett2015}
R.~W. Brockett, {\em Finite Dimensional Linear Systems}.
\newblock Society for Industrial and Applied Mathematics, 2015.

\bibitem{forni2025}
F.~Forni and R.~Sepulchre, ``Gradient modeling of memristive systems,'' in {\em
  Proceedings of the IEEE Conference on Decision and Control}, pp.~25--32,
  IEEE, Dec. 2025.

\end{thebibliography}
\end{document}